\documentclass{article}

\usepackage{amssymb}
\usepackage{latexsym}
\usepackage{mathrsfs}
\usepackage{extarrows}
\usepackage[dvips]{graphicx}

\newtheorem{theo}{Theorem}[section]
\newtheorem{obser}[theo]{Observation}
\newtheorem{remark}[theo]{Remark}

\newtheorem{lemma}[theo]{Lemma}
\newtheorem{claim}[theo]{Claim}
\newtheorem{coro}[theo]{Corollary}
\newtheorem{con}[theo]{Conjecture}
\newtheorem{prop}[theo]{Proposition}
\newtheorem{fact}[theo]{Fact}
\newtheorem{defi}[theo]{Definition}

\def\q{\hspace*{\fill}$\Box$\medskip}

\def\endproofbox{\hskip 1.3em\hfill\rule{6pt}{6pt}}

{%
\noindent{\it Proof}
}%
{%
 \quad\hfill\endproofbox\vspace*{2ex}
}

\begin{document}
\title{Lagrangian densities of hypergraph cycles}
\author{Zilong Yan \thanks{ College of Mathematics and Econometrics, Hunan University, Changsha 410082, P.R. China. Email: zilongyan@hnu.edu.cn.} \and Yuejian Peng \thanks{ Corresponding author. College of Mathematics, Hunan University, Changsha, 410082, P.R. China. Email: ypeng1@hnu.edu.cn. \ Supported in part by National Natural Science Foundation of China (No. 11671124).}
}

\maketitle

\begin{abstract}
The Lagrangian density of an $r$-uniform  hypergraph $F$ is $r!$ multiplying the supremum of the Lagrangians of all $F$-free $r$-uniform  hypergraphs.
For an $r$-graph $H$ with $t$ vertices, it is clear that $\pi_{\lambda}(H)\ge r!\lambda{(K_{t-1}^r)}$. We say that an $r$-unform hypergraph $H$ with $t$ vertices is perfect if $\pi_{\lambda}(H)= r!\lambda{(K_{t-1}^r)}$.
A theorem of Motzkin-Straus implies that all  $2$-uniform graphs are perfect.  It is interesting to explore what kind of hypergraphs are perfect. A hypergraph is linear if any 2 edges have at most 1 vertex in common. We propose the following conjecture: (1) For $r\ge 3$, there exists $n$ such that a linear $r$-unofrm hypergraph with at least $n$ vertices is perfect.
(2) For $r\ge 3$, there exists $n$ such that if $G, H$ are perfect $r$-uniform hypergraphs with at least $n$ vertices, then $G\bigsqcup H$ is perfect.
Regarding this conjecture, we obtain  a partial result:
Let $S_{2,t}=\{123,124,125,126,...,12(t+2)\}$. (An earlier result of Sidorenko states that $S_{2,t}$ is perfect \cite{Sidorenko-89}.) Let $H$ be a perfect  $3$-graph with $s$ vertices. Then $F=S_{2,t}\bigsqcup H$ is perfect if $s\geq 3$ and  $t\geq 3$.

There was no known result on Lagrangian densities of hypergraph cycles and there were 3 unsolved cases for $3$-uniform graphs spanned by 3  edges: a linear cycle of length 3: $C_3^3= \{123, 345, 561\}$, the generalized triangle: $F_5=\{123, 124, 345\}$  and  $K_4^{3-}=\{123, 124, 134\}$. In this paper, we  obtain the Lagrangian density of  $F_5$ and this is the first example of non-perfect $3$-uniform graph.
We also obtain an extension of the above result to $r$-uniform hypergraphs.
We show that $C_3^3$ is perfect,  and among all $C_3^3$-free $3$-graphs $G$, only those hypergraphs containing $K_5^3$ achieves the Lagrangian $\lambda{(K_{5}^3)}$. An extension of this result to the $3$-uniform linear cycle of length $t$ is also given in the paper. The Tur\'an densities of  extensions of the above hypergraphs can be obtained by
applying  a transference technique of Pikhurko.

\end{abstract}


\section{Introduction}
\subsection{ Notations and definitions}
For a set $V$ and a positive integer $r$, $V^ r$ denotes the family of all $r$-subsets of $V$. An {\em $r$-uniform graph} or {\em $r$-graph $G$} consists of a set $V(G)$ of vertices and a set $E(G) \subseteq V(G) ^r$ of edges. Let $|G|$ denote the number of edges of $G$. An edge $e=\{a_1, a_2, \ldots, a_r\}$ will be simply denoted by $a_1a_2 \ldots a_r$. An $r$-graph $H$ is  a {\it subgraph} of an $r$-graph $G$, denoted by $H\subseteq G$, if $V(H)\subseteq V(G)$ and $E(H)\subseteq E(G)$. A subgraph of $G$ {\em induced} by $V'\subseteq V$, denoted as $G[V']$, is the $r$-graph with vertex set $V'$ and edge set $E'=\{e\in E(G):e \subseteq V'\}$. For $S\subseteq V(G)$, let $G\setminus S$ denote the subgraph of $G$ induced by $V(G)\setminus S$.
Let $K^{r}_t$ denote the complete $r$-graph on $t$ vertices. Let $K^{r-}_t$ be the complete $r$-graph on $t$ vertices by removing 1 edge. For a positive integer $n$,  let $[n]$ denote $\{1, 2, 3, \ldots, n\}$.

Given an $r$-graph $F$,  an $r$-graph $G$ is called $F$-free if it does not contain an isomorphic copy of $F$. For a fixed positive integer $n$ and an $r$-graph $F$, the {\em Tur\'an number} of $F$, denoted by $ex(n,F)$, is the maximum number of edges in an $F$-free $r$-graph on $n$ vertices.
 An averaging argument of Katona, Nemetz and Simonovits \cite{KNS} shows that the sequence ${ ex(n,F) \over {n \choose r } }$ is a non-increasing sequence of real numbers in $[0,1]$. Hence, $\lim_{n\rightarrow\infty} { ex(n,F) \over {n \choose r } }$ exists. The {\em Tur\'{a}n density} of $F$ is defined as $$\pi(F)=\lim_{n\rightarrow\infty} { ex(n,F) \over {n \choose r } }.$$

For 2-graphs, Erd\H{o}s-Stone-Simonovits determined the Tur\'an numbers of all non-bipartite graphs asymptotically. Very few results are known for hypergraphs and a recent survey on this topic can be found in Keevash's survey paper \cite{Keevash}.


\begin{defi}
Let $G$ be an $r$-graph on $[n]$ and let
  $\vec{x}=(x_1,\ldots,x_n) \in [0,\infty)^n$. For every subgraph $H\subseteq G$,
define the {\em Lagrangian} function
$$\lambda (H,\vec{x})=\sum_{e \in E(H)}\prod\limits_{i\in e}x_{i}.$$
\end{defi}
The {\em Lagrangian} of
$G$, denoted by $\lambda (G)$, is defined as
 $$\lambda (G) = \max \{\lambda (G, \vec{x}): \vec{x} \in \Delta \},$$
where $$\Delta=\{\vec{x}=(x_1,x_2,\ldots ,x_n) \in [0, 1])^{n}: x_1+x_2+\dots+x_n =1 \}.$$

The value $x_i$ is called the {\em weight} of the vertex $i$ and a vector $\vec{x} \in {\Delta}$ is called a {\em feasible weight vector} on $G$.
A feasible weight vector  $\vec{y}\in {\Delta}$ is called an {\em optimum weight vector} for $G$ if $\lambda (G, \vec{y})=\lambda(G)$.

In \cite{MS}, Motzkin and Straus established a connection between the Lagrangian of any given $2$-graph and it's maximum complete subgraphs.
\begin{theo} {(\cite{MS})} \label{MStheo}
If $G$ is a $2$-graph in which a maximum complete subgraph has  $t$ vertices, then
$\lambda(G)=\lambda(K_t^2)={1 \over 2}(1 - {1 \over t})$.
\end{theo}

 They also applied this connection to give another proof of the theorem of Tur\'an  on the Tur\'an density of complete graphs. Since then, the Lagrangian method has been a useful tool in hypergraph extremal problems. Earlier applications include that Frankl and R\"{o}dl \cite{FR} applied it in disproving the long standing jumping constant conjecture of Erd\H{o}s. Sidorenko \cite{Sidorenko-89}  and Frankl-F\"uredi \cite{FF} applied Lagrangians of hypergraphs in  finding Tur\'an densities of some hypergraphs. More recent developments of this method were obtained by Pikhurko \cite{Pikhurko} and in the papers \cite{HK, NY, BIJ, NY2, Jenssen}.
 In addition to its applications in extremal problems,  it is interesting in its own right to determine the  maximum Lagrangian of $r$-graphs with certain properties as remarked by Hefetz and Keevash \cite{HK}. For example, an interesting conjecture of Frankl-F\"uredi \cite{FF} considers the question of determining the maximum  Lagrangian among all $r$-graphs with the fixed  number of edges.  Talbot \cite{T} made a first breakthrough in  confirming  this conjecture for some cases. Subsequent progress in this conjecture  were made in the papers of Tyomkyn \cite{Tyo}, Lei-Lu-Peng \cite{LLP2018} and Tang-Peng-Zhang-Zhao \cite{TPZZ2}. Recently, Gruslys-Letzter-Morrison \cite{GLM2018} confirmed this conjecture for $r=3$ and sufficiently large $m$.  We focus on  the Lagrangian density of an $r$-graph $F$ in this paper.

Given an $r$-graph $F$, the {\em Lagrangian density } $\pi_{\lambda}(F)$ of $F$ is defined to be
$$\pi_{\lambda}(F)=\sup \{r! \lambda(G): G \:\: {\rm is} \:\: F{\text-}{\rm free}\}.$$

The Lagrangian density  is closely related to the Tur\'an density. It's easy to show the following fact.

\begin{fact}\label{prop2}
$\pi(F)\leq\pi_{\lambda}(F).$
\end{fact}
\noindent{\em  Proof.}
Let $\epsilon>0$ be arbitrary. Let $n$ be large enough and let $G_n$ be a maximum $F$-free $r$-graph on $n$ vertices such that $\pi(F)\leq\frac{|G_n|}{\binom{n}{r}}+\frac{\epsilon}{2}$. Then $$\pi(F)\leq\frac{|G_n|}{\binom{n}{r}}+\frac{\epsilon}{2}\leq r!\sum_{e\in E(G_n)}\frac{1}{n^r}+\epsilon=r!\lambda{(G_n, (\frac{1}{n},\frac{1}{n},\dots,\frac{1}{n}))}+\epsilon\leq r!\lambda{(G_n)}+\epsilon\leq\pi_{\lambda}(F)+\epsilon.$$
\q

A pair of vertices $\{i, j\}$ is {\em covered}  in a hypergraph $F$ if there exists  an edge $e$ in $F$ such that  $\{i, j\}\subseteq e$. We say that $F$ covers pairs if every pair of vertices in $F$ is covered.
 Let $r\ge 3$ and $F$ be an $r$-graph.  The  {\em extension} of $F$, denoted by $H^F$ is
obtained as follows: For each pair of vertices $v_i,v_j$ not covered in $F$, we add a set $B_{ij}$ of $r-2$ new vertices and the edge $\{v_i,v_j\} \cup B_{ij}$, where the $B_{ij}$'s are pairwise disjoint over all such pairs $\{i,j\}$. A transference technique of Sidorenko\cite{Sidorenko-87} and Pikhurko\cite{Pikhurko} gave the following connection between the Lagrangian density of a hypergraph  $F$ and the Tur\'an density of its extension.

\begin{prop}\cite{Sidorenko-87, Pikhurko}\label{prop1}
$\pi(H^F)= \pi_{\lambda}(F).$ In particular, if $F$ covers pairs, then $\pi(F)=\pi_{\lambda}(F).$
\end{prop}

For example, to determine the Tur\'an density of $K_4^{3}$  (a long standing conjecture of Tur\'an) is equivalent to determine the Lagrangian density of $K_4^{3}$.

 The Lagrangian density of the enlargement of a tree satisfying Erd\H{o}s-Sos's conjecture is determined by Sidorenko \cite{Sidorenko-89} and Brandt-Irwin-Jiang \cite{BIJ}.
Pikhurko \cite{Pikhurko} determined the Lagrangian density of a $4$-uniform tight  path of length 2  and this led  to confirm the conjecture of Frankl-F\"uredi on the  Tur\'an number of its extension, the $r$-uniform genearlized triangle for the case $r=4$.   Norin and Yepremyan \cite{NY2} determined  for $r=5$ or $6$  by extending the earlier result of Frankl-F\"uredi in \cite{FF}. Jenssen \cite{Jenssen} determined the Lagrangian density of a path of length 2 formed by two edges intersecting at   $r-2$ vertices for $r=3, 4, 5, 6, 7$.
Hefetz and Keevash \cite{HK} determined the Lagrangian density of a $3$-uniform matching  of size 2. Jiang-Peng-Wu \cite{JPW} obtained for any $3$-uniform matching. The case for  an $r$-uniform matching  of size 2   was given  in \cite{NWY} (independently, in \cite{WPC} for $r=4$). In \cite{WP, HPW, CLP},  we obtained the Lagrangian densities of a $3$-uniform linear path of length $3$ or $4$, the disjoint union of a $3$-uniform linear path of length $2$ or $3$ and a  $3$-uniform matching, and the disjoint union of a $3$-uniform tight path of length $2$  and a  $3$-uniform matching.These were all the known results on Lagrangian densities.

\subsection{ Main results, open problems and Remarks}


 For an $r$-graph $H$ on $t$ vertices, it is clear that $\pi_{\lambda}(H)\ge r!\lambda{(K_{t-1}^r)}$. We say that an $r$-uniform hypergraph $H$ on $t$ vertices is {\em perfect} if $\pi_{\lambda}(H)= r!\lambda{(K_{t-1}^r)}$.
Theorem \ref{MStheo} implies that all  $2$-graphs are perfect.  Currently,  all hypergraphs with known Lagrangian densities are perfect (mentioned in the previous paragraph) except an $r$-uniform matching of size 2 for $r\ge4$. It is interesting to explore what kind of hypergraphs are perfect. An $r$-uniform hypergraph is {\em linear} if any two edges have at most 1 vertex in common. Let $G\bigsqcup H$ denote the disjoint union of $G$ and $H$. Let us  propose the following conjecture.

\begin{con}\label{con1}

(1) For $r\ge 3$, there exists $n$ such that a linear $r$-unofrm hypergraph with at least $n$ vertices is perfect.

(2) For $r\ge 3$, there exists $n$ such that if $G, H$ are perfect $r$-graphs with at least $n$ vertices, then $G\bigsqcup H$ is perfect.
\end{con}

For $r\ge 4$, the condition that the number of vertices is large enough cannot be removed from the above conjecture as the $r$-uniform matching of size 2 is not perfect for $r\ge 4$.

Regarding this conjecture, we provide a partial result.

\begin{theo}\label{theo3}
Let $S_{2,t}=\{123,124,125,126,...,12(t+2)\}$. Let $H$ be a perfect  $3$-graph with $s$ vertices. Then $F=S_{2,t}\bigsqcup H$ is perfect  if $s\geq 3$ and $t\geq 3$.
\end{theo}

An earlier result of Sidorenko states that $S_{2, t}$ is perfect \cite{Sidorenko-89}. Combining previous known results this Theorem generates many perfect hypergraphs.

In view of the known results, there was no  result on hypergraph cycles and there were 3 unsolved cases for $3$-uniform graphs spanned by 3  edges: a linear cycle of length 3 denoted by $C_3^3= \{123, 345, 561\}$, the generalized triangle denoted by $F_5=\{123, 124, 345\}$ (a Berge cycle) and  $K_4^{3-}=\{123, 124, 134\}$ (a tight cycle). In this paper, we also obtain the Lagrangian density of  $F_5$ and $C_3^3$. We  prove the following result implying that $F_5$ is not perfect.

 \begin{theo}  \label{3.4}\label{theo1}
Let $G$ be an $F_5$-free 3-graph with $n$ vertices.   Then $\lambda{(G)}\leq\frac{2}{27}$. Furthermore, equality holds if and only if $G$ contains a copy of $S_{n'}^3(1)$ with $n'\leq n$ and $n'\to\infty$, where $S_n^3(1)$ is the $3$-graph with vertex set $[n]$ and edge set $\{1ij| \{i, j\}\subset [n]\setminus \{1\}\}$.
\end{theo}

 \begin{coro}  \label{coro1}
$\pi_{\lambda}(F_5)=\frac{4}{9}.$
\end{coro}
\noindent{\em  Proof of Corollary \ref{coro1}.}
Since $S_n^3(1)$ is $F_5$-free, then $\pi_{\lambda}(F_5)\geq\lim\limits_{n\to\infty}3!\lambda{(S_n^3(1))}=\frac{4}{9}.$ On the other hand, by Theorem \ref{3.4}, $\pi_{\lambda}(F_5)\leq\frac{4}{9}.$ Therefore, $\pi_{\lambda}(F_5)=\frac{4}{9}.$
\q

We also consider an extension of the above result to $r$-graphs.
Let edges $e_1=\{1,2,3,...,r-2,a_1,a_2\}$,  $e_2=\{1,2,3,\dots,r-2,a_1,a_3\}$,\dots,  $e_{s-1}=\{1,2,3,\dots,r-2,a_1,a_s\}$. 
Let $r$-uniform graphs $F_0^{r}=e_1\cup e_2\cup \{a_2,a_3,m_1,m_2, \dots,m_{r-2}\}$, $F_i^{r}=e_1\cup e_2\cup \{a_2,a_3,1,2,\dots i,m_1,m_2,\dots,m_{r-i-2}\}$ for $1\leq i\leq r-3$. Let $ \mathscr{F}^{r}=\{F_0^{r}, F_1^{r},\dots,F_{r-3}^{r}\}$. Note that  $F_0^{3}$ is $F_5.$ We show that
\begin{theo} \label{theo1'}
Let $s=\lceil \frac{r^{r-1}}{2(r-1)!}\rceil$ and  $G$ be an $\mathscr{F}^{r}$-free $r$-graph satisfying $e_1,e_2,\dots,e_{s-1}\in E(G)$. Then $\lambda{(G)}\leq\frac{2}{r^{r}}$.
\end{theo}

As remarked in Remark \ref{3.5}, when $r=3$, Theorem \ref{theo1'} implies the first part of Theorem \ref{theo1}.

An $r$-uniform linear cycle of length $t$ denoted by $C_t^3$ is isomorphic to $\{123, 345,\dots, (2t-1)(2t)1\}$.  Since $C_t^3$ has $2t$ vertices, then $K_{2t-1}^3$ is $C_t^3$-free, consequently, $\pi_{\lambda}(C_t^3)\geq 3!\lambda{(K_{2t-1}^3)}$. We  show that $C_3^3$ is perfect, and among all $C_3^3$-free $3$-graphs $G$, only those containing $K_5^3$ achieves the Lagrangian $\lambda{(K_{5}^3)}$.
\begin{theo}\label{theo2}\label{c33}
Let $G$ be a $C_3^3$-free $3$-graph. Then $\lambda{(G)}\leq\lambda{(K_5^3)}=\frac{2}{25}$, and equality holds if and only if $G$ contains $K_5^3$ as a subgraph.
\end{theo}

\begin{coro}  \label{coro2}
$\pi_{\lambda}(C_3^3)=\frac{12}{25}.$
\end{coro}
\noindent{\em  Proof of Corollary \ref{coro2}.}
Since $K_5^3$ is $C_3^3$-free and $\lambda{(K_5^3)}=\frac{2}{25}$, then $\pi_{\lambda}(C_3^3)\geq3!\lambda{(K_5^3)}=\frac{12}{25}.$ On the other hand, by Theorem \ref{theo2}, $\pi_{\lambda}(C_3^3)\leq\frac{12}{25}.$ Therefore, $\pi_{\lambda}(C_3^3)=\frac{12}{25}.$
\q

An $r$-graph $G$ is {\em dense} if $\lambda (G') < \lambda (G)$ for every proper subgraph $G'$ of $G$. As remarked in Remark \ref{remark}, if $G$ is $F$-free r-graph, then $G$ contains a dense subgraph with the same Lagrangian as $G$. So to estimate an upper bound of the Lagrangians of all $C_t^{3}$-free 3-graphs, we only need to consider all dense $C_t^{3}$-free 3-graphs. The following result is an extension of Theorem \ref{c33} (this will be explained in Section \ref{sec5}).

\begin{theo}\label{3t'}
Let $G$ be a dense and $C_t^{3}$-free 3-graph satisfying $K_{2t-2}^{3-}\subseteq G$ $(t\geq 4)$, then $\lambda{(G)}\leq \lambda{(K_{2t-1}^3)}.$
\end{theo}

The condition that $K_{2t-2}^{3-}\subseteq G$ in Theorem \ref{3t'} is very strong, with effort, this condition can be weakened (We omit the technical argument. Instead, only the proof of Theorem \ref{3t'} will be given).  We propose the following conjecture.

\begin{con}
$C_t^{3}$ is perfect for any $t\ge 3$. A $C_t^{3}$-free $3$-graph achieves the maximum Lagrangian density only on $3$-graphs containing $K_{2t-1}^3$.
\end{con}

The first part of the conjecture is included in Conjecture \ref{con1}.

In Section \ref{sec2}, we give some preliminary results on the Lagrangian function. In Section \ref{sec3}, we give the proofs of Theorems \ref{theo1} and \ref{theo1'}. The proof of Theorem \ref{theo3} will be given in Section \ref{sec4}, and the proof of Theorem \ref{c33} and \ref{3t'} will be given in Section \ref{sec5}.

{\bf Remark.} The Tur\'an densities of  extensions of the above hypergraphs can be obtained by
applying Proposition \ref{prop1}. 
There is still one unsolved case for $3$-uniform graphs spanned by 3  edges:  $K_4^{3-}=\{123, 124, 134\}$. Since $K_4^{3-}$ covers pairs, the Lagrangian density of $K_4^{3-}$ is the same as the  Tur\'an density of $K_4^{3-}$. It would be very interesting if the Tur\'an density of $K_4^{3-}$ can be obtained by determining the Lagrangian density of $K_4^{3-}$.

\section{Preliminaries}\label{sec2}

The following fact follows immediately from the definition of the Lagrangian.
\begin{fact}\label{mono}
Let $G_1$, $G_2$ be $r$-graphs and $G_1\subseteq G_2$. Then $\lambda (G_1) \le \lambda (G_2).$
\end{fact}

\begin{fact} {\em (\cite{FR})}\label{fact2}
Let $G$ be an $r$-graph on $[n]$. Let $\vec{x}=(x_1,x_2,\dots,x_n)$ be an optimum weight vector on  $G$. Then
$$ \frac{\partial \lambda (G, \vec{x})}{\partial x_i}=r\lambda(G)$$
for every $i \in [n]$ satisfying $x_i>0$.
\end{fact}

This can be generalized to the following result.
\begin{fact}
Let $E\subset [n]^r$ and $f(\overrightarrow{x})$ be a homogeneous function with degree $r$ in the form of $\sum\limits_{\{i_1,i_2,\dots,i_r\}\in E}$ $a_{i_1i_2\dots i_r}x_{i_1}x_{i_2}\dots x_{i_r}$. Let $\vec{x}=(x_1,x_2,\dots,x_n)$ be an optimum weight vector for  $Max \{f(\overrightarrow{x}), \overrightarrow{x}\in S\}$. Then
$$ \frac{\partial f(\overrightarrow{x})}{\partial x_i}=rf(\overrightarrow{x})$$
for every $i \in [n]$ satisfying $x_i>0$.
\end{fact}

Given an $r$-graph $G$, and $i, j\in V(G),$ define $$L_G(j\setminus i)=\{e: i\notin e, e\cup\{j\}\in E(G)\:and\: e\cup\{i\}\notin E(G)\}.$$

\begin{fact}\label{symmetry}
Let $G$ be an $r$-graph on $[n]$. Let $\vec{x}=(x_1,x_2,\dots,x_n)$ be a feasible weight vector on $G$. Let $i,j\in [n]$, $i\neq j$ satisfying $L_G(i \setminus j)=L_G(j \setminus i)=\emptyset$. Let
$\vec{y}=(y_1,y_2,\dots,y_n)$ be defined by letting $y_\ell=x_\ell$ for every $\ell \in [n]\setminus \{i,j\}$ and $y_i=y_j={1 \over 2}(x_i+x_j)$.
Then $\lambda(G,\vec{y})\geq \lambda(G,\vec{x})$. Furthermore, if the pair $\{i,j\}$ is covered by an edge of $G$, $x_i>0$ for each $1\le i\le n$,  and $\lambda(G,\vec{y})=\lambda(G,\vec{x})$, then $x_i=x_j$.
\end{fact}
\noindent{\em  Proof.}
Since $L_G(i \setminus j)=L_G(j \setminus i)=\emptyset$, then
$$\lambda(G,\vec{y})-\lambda(G,\vec{x})=\sum_{\{i,j\} \subseteq e \in G}\left({(x_i+x_j)^2 \over 4}-x_ix_j\right)\prod\limits_{k\in e\setminus \{i,j\}}x_k \ge 0.$$
If the pair $\{i,j\}$ is covered by an edge of $G$ and $x_i>0$ for each $1\le i\le n$, then the equality holds only if $x_i=x_j$.
\q


\begin{fact} {\em (\cite{FR})}\label{dense}
Let $G=(V,E)$ be a dense $r$-graph. Then $G$ covers pairs.
\end{fact}

Note that the converse of Fact \ref{dense} is not true.  For example, the Fano plane covers pairs but it is not dense. Indeed, many counterexamples exist by Theorem 2.1 in the paper of Talbot \cite{T}.
\begin{fact} \label{2.6} {\em (\cite{Sidorenko-87})}
If $G$ is a dense 3-graph on $[n]$ $(n\geq 4)$. Then $G$ contains a copy  isomorphic to $\{123, 124\}$.
\end{fact}

While considering the Lagrangian density of an $r$-graph $F$, we can always reduce to consider dense $F$-free $r$-graph.
\begin{remark} \label{remark}
Let $F$, $G$ be $r$-graphs and $G$ be $F$-free. Then there exists a dense subgraph $G'$ of $G$ such that $\lambda{(G')}=\lambda{(G)}$ and $G'$ is $F$-free.
\end{remark}
\noindent{\em  Proof.}
Let $G$ be an $r$-graph on $n$ vertices. If $G$ is dense, then we are fine. If not, then we can find $G'\subset G$ such that $\lambda{(G')}=\lambda{(G)}$ and $|V(G')|<|V(G)|.$ If $G'$ is dense, then we stop. Otherwise, we continue this process until we find a dense subgraph.
\q


\section{Lagrangian density of generalized triangle}\label{sec3}

In an $r$-graph $G$, $N(a)$ denotes the link of $a$, i.e. $N(a)=\{S|\{a\}\cup S\in E(G)\}$. Let $N(a,b)$ denote the link of $\{a, b\}$, i.e. all $(r-2)$-sets $S$ such that $\{a, b\}\cup S\in E(G)$. Let $N_A(a)$ denote $N(a)\cap A$.

\subsection{Lagrangian density of $F_5$}

 In the section, the proof of Theorem \ref{3.4} will be given. Let $G$ be an $F_5$-free $3$-graph. By Remark \ref{remark}, we may assume that $G$ is dense, $F_5$-free, and $|G|>5$.  We first show several crucial facts.
\begin{fact}\label{compression-dense}
Let $G$ be a dense graph on n vertices and let $\vec{x}=(x_1, x_2, ..., x_n)$ be an optimal weight vector satisfying $x_i\geq{x_{i+1}}>0$. If there exist $x_i$ such that $x_i\geq\frac{1}{3}$, then\\
(i) $\lambda(G)\leq\frac{2}{27}$.\\
(ii) $\lambda(G)\rightarrow\frac{2}{27}$ if and only if $x_1=\frac{1}{3}$ and $G\supseteq S_t^3(1)$ satisfying $t\to\infty$.
\end{fact}
\noindent{\em  Proof.}
Since $x_1\geq{x_i\geq\frac{1}{3}}$, then $(1-x_1)^2\leq\frac{4}{9}$. By Theorem \ref{MStheo} and Fact \ref{fact2}, $$3\lambda{(G)}=\frac{\partial \lambda{(G)}}{\partial {x_1}}\leq{(\frac{1-x_1}{t-1})^2\binom{t-1}{2}}=\frac{(t-2)(1-x_1)^2}{2(t-1)}\leq\frac{2}{9},$$ where $t$ is the maximum complete subgraph in $N(1).$

So $\lambda(G)\leq\frac{2}{27}$. And $\lambda(G)\rightarrow\frac{2}{27}$ if and only if $x_1=\frac{1}{3}$, $G\supseteq S_t^3(1)$ satisfying $t\to\infty$.
\q

\begin{lemma}\label{lemma1}
If $G$ is dense, $F_5$-free and contains a $K_4^3$, then $\lambda{(G)}\leq\frac{1}{16}$.
\end{lemma}
\noindent{\em  Proof.}
Assume that $v_1,v_2,v_3,v_4\in V(G)$ form a $K_4^3$. We first prove
\begin{claim} \label{claim}
Let $v_i, v_j\in V(G)\setminus \{v_1, v_2, v_3, v_4\}$. Let $k, t\in [4]$, and $k\neq t$, then the following properties hold:\\
  (i) if $v_iv_jv_k\in E(G)$, then $v_iv_jv_t\notin E(G)$,\\
  (ii) $v_iv_kv_t\notin E(G)$.
\end{claim}
\noindent{\em  Proof.}
(i) If $v_iv_jv_k, v_iv_jv_t\in E(G)$,  let $r\in [4]\setminus \{k,t\}$,  then $v_iv_jv_k, v_iv_jv_t, v_kv_tv_r$ form a $F_5$, a contradiction.
\\(ii) If $v_iv_kv_t\in E(G)$, let $p,r\in [4]\setminus \{k, t\}$, then $v_iv_kv_t, v_pv_rv_k, v_pv_rv_t$ form a $F_5$, a contradiction.
\q
\\Let us continue the proof of Lemma \ref{compression-dense}. Let $\vec{x}$ be an optimal weight vector for $G$, and $x_1,x_1,x_3,x_4$ be the weights of $v_1,v_2,v_3,v_4$ respectively. Let $a=x_1+x_2+x_3+x_4$, then by Fact \ref{fact2},
$$12\lambda{(G)}=\sum_{i=1}^4\frac{\partial \lambda{(G)}}{\partial {x_i}}=\sum_{i=1}^4\lambda{(N(v_i), \overrightarrow{x})}.$$
Since $\sum\limits_{i=1}^4\lambda{(N(v_i), \overrightarrow{x})}$ contains each $x_ix_j$ $(5\leq i<j\leq n)$ at most once by Claim \ref{claim}(i), $x_ix_j$ $(1\leq i<j\leq 4)$ exactly twice, and no term in the form of $x_ix_j$ $(1\leq i\leq 4, 5\leq j\leq n)$ by Claim \ref{claim}(ii). Then
$$12\lambda{(G)}\leq\bigg(\frac{1-a}{n-4}\bigg)^2\binom{n-4}{2}+2(x_1x_2+x_1x_3+x_1x_4+x_2x_3+x_2x_4+x_3x_4).$$
Since $\sum\limits_{1\leq i<j\leq 4}x_ix_j\leq6(\frac{a}{4})^2=\frac{3}{8}a^2$, then $$\lambda{(G)}\leq\frac{(1-a)^2}{24}+\frac{a^2}{16}=\frac{5a^2-4a+2}{48}.$$ The above quadratic function obtains the maximum $\frac{1}{16}$ at $a=1$ (recall that $a\in [0,1]$). So $\lambda{(G)}\leq\frac{1}{16}$.
\q

\noindent{\em  Proof of Theorem \ref{3.4}.}
By Remark \ref{remark} and Lemma \ref{lemma1} we may assume that $G$ is a dense $K_4^3$-free 3-graph with $n\geq 5$. By Fact \ref{2.6}, we may assume that $G$ contains $\{012,013\}$. Note that $N(2,3)\subset \{0,1\}$ since otherwise we will get a $F_5$. Without loss of generality, let $023 \in E(G)$. Note that $\{012, 013, 023\}$ form a $K_4^{3-}$.

Let $A$ a maximal set containing $\{1,2,3\}$ such that $0ij\in E(G)$ for any pair $\{i,j\}$ in $A$. Clearly $|A|\geq3$. We divide $V(G)$ into 3 parts, $V(G)=\{0\}\cup{A}\cup{B}$, where $B=V(G)\setminus A\cup\{0\}$. Then we have the following Observation.

\begin{obser}\label{obser1}\label{obser2}\label{obser3}\label{obser4}\label{obser5}
(i) For $v\in B$ and $i, j\in A$, $vij\notin E(G)$.\\
(ii) Let $v\in B$, if $v0i\in E(G)$ for some $i\in A$, then $v0j\in E(G)$ for each $j\in A$. \\
(iii) For any $v\in B$ and any $i\in A$, $\{v, i\}$ is not covered in $G$\\
(iv) For any $v_1,v_2\in B$ we have $|\{i|iv_1v_2\in E(G),i\in A\bigcup\{0\}\}|\leq1$.\\
(v) If $i,j,k\in A$, then $ijk\notin E(G)$.
\end{obser}
\noindent{\em  Proof.}
(i) If $vij\in E(G)$ then $vij$, $0ki$, $0kj$ form a $F_5$ for any $k\in A\setminus \{i,j\}$, such a $k$ exists since $|A|\geq 3.$

(ii) Since $G$ is dense then $N(v, j)\neq \emptyset $ for each $j\in A$. If $N(v, j)\bigcap B\neq \emptyset $, let $u\in N(v, j)\bigcap B$, then $uvj, v0i, 0ij$ form a $F_5$,  so $N(v, j)\bigcap B=\emptyset $. By (i), $N(v, j)\cap A=\emptyset$. So $N(v, j)=\{0\}$.

(iii) This is due to (i), (ii) and the maximality of $A$.

(iv) If there exist $i,j\in A\bigcup\{0\}$ s.t. $iv_1v_2,jv_1v_2\in E(G)$, then $iv_1v_2,jv_1v_2,0ij$ form a $F_5$ (change $0ij$ to $ijk$ when $i=0$ or $j=0$ for $k\in A\setminus \{i,j\}$).

(v) If $ijk\in E(G)$ then $0ij$, $0ik$, $0jk$, $ijk$ form a $K_4^3$ , a contradiction.
\q




Let us continue the proof of Theorem \ref{3.4}. Let $k=|A|$, the weight of 0 be $c$, the weight of vertices in $A$ be $a_1,a_2,...,a_k$ in the optimal weighting vector. Let $a=a_1+a_2+...+a_k$, consider $\frac{\partial\lambda{(G)}}{\partial x_i}$, $i\in A$ and $\frac{\partial\lambda{(G)}}{\partial x_0}$. By Facy \ref{fact2},
$$(3k+3)\lambda{(G)}=\sum_{i=1}^k\frac{\partial\lambda{(G)}}{\partial x_i}+\frac{\partial\lambda{(G)}}{\partial x_0}.$$
By Observation \ref{obser4}(iv), the terms in the form of $x_ix_j$ where $i, j\in B$ appear at most once in $\sum\limits_{i=1}^k\frac{\partial\lambda{(G)}}{\partial x_i}+\frac{\partial\lambda{(G)}}{\partial x_0}$, the terms in the form of $a_ia_j (1\leq i<j\leq k)$ appear exactly once in $\frac{\partial\lambda{(G)}}{\partial x_0}$ and not in $\sum\limits_{i=1}^k\frac{\partial\lambda{(G)}}{\partial x_i}$ by the definition of $A$ and Observation \ref{obser5}(iii)(v), the terms in the form of $a_ic$ $(1\leq i\leq k)$ appear exactly $k-1$ times in $\sum\limits_{i=1}^k\frac{\partial\lambda{(G)}}{\partial x_i}$ by the definition of $A$. By Observation \ref{obser3}(iii), no term in the form of $cx_j$ $(j\in B)$ or $a_ix_j$ $(1\leq i\leq k, j\in B)$ appear in $\sum\limits_{i=1}^k\frac{\partial\lambda{(G)}}{\partial x_i}+\frac{\partial\lambda{(G)}}{\partial x_0}.$  So
\begin{eqnarray}
(3k+3)\lambda{(G)}&\leq&\binom{n-k-1}{2}\bigg(\frac{1-c-a}{n-k-1}\bigg)^2+\sum\limits_{1\leq i<j\leq k} a_ia_j+(k-1)ac\nonumber\\
&\leq&\frac{(1-a-c)^2}{2}+\frac{1}{2}a^2+(k-1)ac.\nonumber
\end{eqnarray}

So
\begin{eqnarray}
\lambda{(G)}&\leq&\frac{1}{6(k+1)}(c^2+a^2+1+2ac-2c-2a)+\frac{k-1}{3(k+1)}ac+\frac{1}{6(k+1)}a^2\nonumber\\
&\leq&\frac{c^2+2a^2+1+2kac-2c-2a}{6(k+1)}\nonumber\\&=&f(c)\quad(c\leq\frac{1}{3}) .\nonumber
\end{eqnarray}

Since $f(c)$ is quadratic, opens up and $0\leqslant c \leqslant\frac{1}{3}$,  then
\begin{center}
$\lambda{(G)}\leq max\{f(0); f(\frac{1}{3})\}=max\{\frac{2a^2-2a+1}{6(k+1)}(a<1); \frac{2a^2+(\frac{2k}{3}-2)a+\frac{4}{9}}{6(k+1)}  (a\leqslant\frac{2}{3})\}.$
\end{center}

Thus $\lambda{(G)}\leq max\{\frac{1}{6(k+1)},\quad \frac{2}{27}\frac{k}{k+1}\}$.  Since $k\geq 3$, then $\lambda{(G)}\leq\frac{2}{27}$. The proof implies that the  equality holds if and only if $c=\frac{1}{3}$ and $a=\frac{2}{3}$. This implies that the weight of each vertex in $B$ is 0. Then $G$ is dense implies that $B=\emptyset$ and $G$ is isomorphic to $S_n^3(1)$ with $n\to\infty$.
\q

\bigskip
\subsection{Extension to $r$-uniform graphs}

The proof of Theorem \ref{theo1'} will be given in this section.
\begin{remark}\label{3.5}
When $r=3$, Theorem \ref{theo1'} is equivalence to the first part of Theorem \ref{theo1}.
\end{remark}
\noindent{\em  Proof.}
For $r=3$, $s=\frac{3^2}{2\times2!}=\frac{9}{4}\leq 3$. By Fact \ref{2.6} and Remark \ref{remark}, $G$ contains a copy of $\{123, 124\}$.
\q

\noindent{\em  Proof of Theorem \ref{theo1'}.}
By Remark \ref{remark} we may assume that G is dense. Apply induction on $r$. By Theorem \ref{theo1}, the conclusion holds for $r=3$. Assume that the conclusion holds for all the integers less than $r(>3)$. Since $G$ is dense, then there exist an edge $e'\in E(G)$ which covers $\{a_2,a_3\}$. Note that $e'\subset V(e_1\cup e_2)$ since otherwise $e_1\cup e_2\cup e'\in \mathscr{F}^{r}$. Without loss of generality, let $e'=\{1,2,3,...,r-2,a_2,a_3\}$. Similarly there exist an edge $e''$ covering $a_3,a_4$. We note that $e''\subset V(e_2\cup e_3)$ and $e''=\{1,2,3,...,r-2,a_3,a_4\}$ since otherwise $e_1\cup e'\cup e''\in \mathscr{F}^{r}$. So we can divide $V(G)$ into 3 parts. Let $O=\{1,2,3,...,r-2\}$, $\{a_1,a_2,..., a_s\}\subset A$ and $A=\{a_1,a_2,..., a_k\}$ be a maximal set such that for any $a_i,a_j \in A$, $123...(r-2)a_ia_j\in E(G)$. Let $B=V(G)\setminus (O\cup A)=\{b_1,b_2,...\}$, then we have the following facts.
\begin{fact} \label{a}
$O^i\times A^j\times B^l\cap E(G)=\emptyset$, where $j\geq 2, l\geq 1, i+j+l=r$, and $O^i\times A^i\times B^l$ denotes the set of all r-sets consisting of $i$ elements from  $O$, $j$ elements from $A$ and $l$ elements from $B$.
\end{fact}
\noindent{\em  Proof.}
If 
 $e=123...ia_1...a_jb_1...b_l\in O^i\times A^j\times B^l\cap E(G)$, then $\{1,2,...,r-2,a_2,a_3\}\cup \{1,2,...,r-2,a_1,a_3\}\cup e\in \mathscr{F}^{r}$.
\q
\begin{fact}
$O^{r-2}\times A\times B\bigcap E(G)=\emptyset$.
\end{fact}
\noindent{\em  Proof.}
Let $e\in O^{r-2}\times A\times B\bigcap E(G)$. For each $i$, consider $e$ and $e_{i-1}$, then there exist an edge $e'$ covering $\{b_1,a_i\}$, and $V(e')\subset V(e\cup e_i)$. If $e'\neq12...r-2a_ib_1$, then $a_1\in e'$, and $e'\cup$ $\{1, 2, 3, ... , a_i, a_j\}\cup$ $\{1, 2, ..., r-2, a_1, a_j\}\in \mathscr{F}^{r}$ for $j\in [k]\setminus \{i\}$. So $123...(r-2)a_ib_1\in E(G)$ for each $i\in [k]$, then $b_1\in A$ by the maximality of $A$. A contradiction.
\q

Let $\{x_1, x_2, ..., x_{r-2}, x_{a_1}, ... ,x_{a_k}, x_{b_1}, ... ,x_{b_i},...\}$ be an optimal weight of $V(G)=O\cup A\cup B$. If there exist $i\in V(G)$ such that $x_i\geq \frac{1}{r}$, then $N(i)$ is an $\mathscr{F}^{r-1}$-free $(r-1)$-graph. By Fact \ref{fact2} and the induction hypothesis, we have $$r\lambda{(G)}=\frac{\partial\lambda{(G)}}{\partial x_i}\leq\frac{2}{(r-1)^{(r-1)}}\bigg(1-\frac{1}{r}\bigg)^{r-1}=\frac{2}{r^{r-1}}.$$
Consequently $\lambda{(G)}\leq \frac{2}{r^r}$.

So we may assume that $0<x_i<\frac{1}{r}$ for each $i\in V(G)$. By Fact \ref{fact2} we have $$kr\lambda{(G)}=\sum_{i=1}^k \frac{\partial\lambda{(G)}}{\partial x_{a_i}}.$$

Now let function $g(\overrightarrow{x})=\sum\limits_{i=1}^k\frac{\partial\lambda{(G)}}{\partial x_{a_i}}$. Consider the optimization problem $$Max\{g(\overrightarrow{x})|\sum_{i\in V(G)}x_i=1, 0\leq x_i\leq \frac{1}{r}\},$$
and let $\overrightarrow{x'}$ be a solution to the optimization problem. Then we know that $$kr\lambda{(G)}\leq g(\overrightarrow{x'}).\eqno (3.1)$$

If there exist $a_j\in A$ such that $x'_{a_j}>0$, then by Fact \ref{fact2}, $$(r-1)g(\overrightarrow{x'})=\frac{\partial g(x)}{\partial x_{a_j}}\bigg|_{\overrightarrow{x'}}=(k-1)x'_1x'_2...x'_{r-2}.$$\\So $$kr(r-1)\lambda{(G)}\leq (k-1)x'_1x'_2...x'_{r-2}\leq (k-1)\bigg(\frac{1}{r}\bigg)^{r-2}.$$
Thus $$\lambda{(G)}\leq \frac{2}{r^r}.$$

So we can assume that $x'_{a_1}=x'_{a_2}=...=x'_{a_k}=0$. Let us estimate $g(\overrightarrow{x'})$. In this case, all terms containing $x'_{a_i}$ for some $a_i\in A$ are 0. We only need to look at terms corresponding to $(r-1)$-tuples in $O\cup B$. Note that $N_{O\cup B}(a_i)\cap N_{O\cup B}(a_j)=\emptyset$ since if $12...tb_1b_2...b_{r-t-1}a_i$, $12...tb_1b_2...b_{r-t-1}a_j\in E(G)$, then these 2 edges together with $12...(r-2)a_ia_j$ form an $r$-graph in $\mathscr{F}^{r}$. So the set of subscripts of all non-zero terms in $g(\overrightarrow{x'})$ is a subset of the complete $(r-1)$-graph on $O\cup B$. Therefore, $$g(\overrightarrow{x'})\leq \bigg(\frac{1}{n-k}\bigg)^{r-1}\binom{n-k}{r-1}\leq \frac{1}{(r-1)!}.\eqno(3.2)$$
Combining (3.1) and (3.2), we have  $$\lambda{(G)}\leq \frac{1}{kr}g(\overrightarrow{x'})\leq \frac{1}{kr(r-1)!}\leq \frac{2}{r^r}$$
since $k\geq s\geq \frac{r^{r-1}}{2(r-1)!}.$
\q

\section{Generating perfect hypergraphs}\label{sec4}

The  proof of Theorem \ref{theo3} will be given in this section.

\bigskip

\noindent{\em  Proof of Theorem \ref{theo3}.}
Let $H$ be a perfect 3-graph with $s$ vertices. Let $F=S_{2,t}\sqcup H$. Let $G$ be an $F$-free 3-graph such that $\lambda{(G)}\rightarrow Sup\{\lambda{(G')}|$ $G'$ is F-free$\}$. By Remark \ref{remark}, we may assume that $G$ is dense. If $F$ is not perfect, then $\lambda{(G)}>\lambda{(K_{s+t+1}^3)}$ since the number of vertices in $F$ is $s+t+2$. By a result of Sidorenko (\cite{Sidorenko-89}), $S_{2,s+t}$ is perfect. So $S_{2,s+t}\subseteq G$. Assume that $\{12a_1,12a_2,...12a_{s+t}\}\subseteq G$. We show the following claim.
\begin{claim}
$G\setminus\{1,2\}$ is $H$-free.
\end{claim}
\noindent{\em  Proof.}
If $H\subset G\setminus\{1,2\}$, we know that $S_{2,s+t}\subseteq G$, so there exist at least $t$ vertices $\{a_{i_1}, a_{i_2}, ... , a_{i_t}\}$ in $G\setminus(\{1,2\}\bigcup H)$ such that $12a_{i_1}, 12a_{i_2}, ...,a_{i_t}\in E(G)$ since $|V(H)|=s$. So $G$ is not $F$-free, a contradiction.
\q

Let us continue the  proof of Theorem \ref{theo3}.  Since $G\setminus \{1,2\}$ is $H$-free and $H$ is perfect, then $\lambda{(G\setminus \{1,2\})}\leq \lambda{(K_{s-1}^3)}$. Let $G'$ be the 3-graph with edge set $E(G)\cup \{12i: i\in V(G)\setminus\{1, 2\}\}\cup \{1jk, 2jk: j, k\in V(G)\setminus\{1, 2\}\}$. Clearly $G\subseteq G'.$ By Fact \ref{mono}, we know that $\lambda{(G)}\leq \lambda{(G')}$. By Fact \ref{symmetry}, we can assume that, the vertices 1, 2 have the same weights $a$ in an optimal weight vector of $G'$. Then we have,
\begin{eqnarray}
\lambda{(G)}&\leq&\lambda{(G')}\nonumber\\
&\leq&a^2(1-2a)+2a\bigg(\frac{1-2a}{n-2}\bigg)^2\binom{n-2}{2}+\lambda{(K_{s-1}^{3})}(1-2a)^3\nonumber\\
&\leq&2a^3-3a^2+a+\frac{(1-2a)^3}{6}m=f(a)\quad (m=\frac{(s-2)(s-3)}{(s-1)^2})\nonumber.
\end{eqnarray}
By direct calculation $f'(a)=(6-4m)a^2+(4m-6)a+1-m$ and
$f''(a)=2(6-4m)(a-\frac{1}{2})\leq0 $ when $a\leq\frac{1}{2}.$
So $f'(a)$ is monotone decreasing in $(0, \frac{1}{2})$. Since $f'(\frac{1}{2}\pm \frac{\sqrt{3-2m}}{6-4m})=0$, then $f(a)$ is monotone increasing in $(0,\frac{1}{2}-\frac{\sqrt{3-2m}}{6-4m})$, monotone decreasing in $(\frac{1}{2}-\frac{\sqrt{3-2m}}{6-4m}, \frac{1}{2})$. Since $a\in (0, \frac{1}{2})$, then  $$f(a)\leq f(\frac{1}{2}-\frac{\sqrt{3-2m}}{6-4m})=\frac{1}{6\sqrt{3-2m}}=\frac{s-1}{6\sqrt{s^2+4s-9}}.$$ So $$\lambda{(G)}<\frac{s-1}{6\sqrt{s^2+4s-9}}.$$

Our goal is to show that $\lambda{(G)}<\lambda{(K_{s+t+1}^{3})}.$ Since for $t\geq 3$
$$\frac{(s+3)(s+2)}{6(s+4)^2}=\lambda{(K_{s+4}^{3})}\leq\lambda{(K_{s+t+1}^{3})}.$$ So it's sufficient to show that $\frac{s-1}{6\sqrt{s^2+4s-9}}\leq\frac{(s+3)(s+2)}{6(s+4)^2}.$ This is equivalent to show that
$$(s-1)(s+4)^2<(s+3)(s+2)\sqrt{s^2+4s-9}.$$ i.e. $$(s^2-2s+1)(s^2+8s+16)^2<(s^2+5s+6)^2(s^2+4s-9).$$
Subtracting the left hand side in both sides, we got
$$0<3s^4+38s^3+103s^2-140s-580=g(s).$$ So our goal is to show the above inequality $g(s)>0.$  Since $g'(s)=12s^3+114s^2+206s-140$, and $g''(s)=36s^2+228s+206 \quad (s>0),$
then $g''(s)>0$ for all $s>0$. So $g'(s)$ is monotone increasing. By direct calculation  $g'(0)<0$, $g'(1)>0$, so $g(s)$ is monotone increasing in $[1,+\infty)$. Since $g(3)>0$, then $g(s)>0$ holds for all $s\geq 3$. This completes the proof.
\q


\section{Lagrangian density of a 3-uniform linear cycle}\label{sec5}
Let $P_t^3$ denote the $3$-uniform linear path with length $t$.  Let $G$ be a $C_t^3$-free $3$-graph and $A\subset V(G)$. By Remark \ref{remark}, we may assume that $G$ is dense. We will give some structure analysis on $G$  and show that $G$ is a `good' $3$-graph with some nice structure.

\begin{defi}
For two vertices $a_i, b_i\in V(G)\setminus A$, we say that $\{a_i, b_i\}$ is a good pair to $A$ if $N(a_i, k)=b_i$ and $N(b_i, k)=a_i$ for all $k\in A$.
\end{defi}

\begin{defi}
We say that a dense 3-graph $G$ is a good graph to $A$ if $V(G)=\cup_{i=1}^s\{a_i,b_i\}\cup A$ and $\{a_i,b_i\}$ is a good pair to $A$ for each $1\leq i \leq s$, where $s=\frac{|V(G)\setminus A|}{2}$.
\end{defi}

\begin{fact} \label{3t1}
Let $G$ be a dense and $C_t^{3}$-free graph. Let $[2t-2]^3\setminus \{123\}$ be a $K_{2t-2}^{3-}$ in $G$.
Let $a, b, c\in V(G)\setminus [2t-2]$ and $i, j, k, l\in [2t-2]$ be four different integers. Then the following properties hold.\\
(i) If $t\geq4$, then there is at most one of $aij$ and $bkl$ in $E(G)$.\\
(ii) There is at most one of $abi$ and $acj$ in $E(G)$,\\
(iii) There is at most one of $aij, abk$ in $E(G)$.
\end{fact}
\noindent{\em  Proof.}
(i)  If there exist $i, j, k, l\in [2t-2]$ such that $aij\in E(G)$ and $bkl\in E(G)$, then there exist $u\in [2t-2]$, such that $iku\in E(G)$ and $\{i,k,u\}\bigcap \{1,2,3\}\neq \emptyset$, then there exists a linear path $P_{t-3}$ in $[2t-2]\setminus\{i, k, u\}$ with endpoints $j, l$ for $t\ge 4$. Then $P_{t-3}, aij, iku, bkl$ form a $C_{t}^{3}$, a contradiction.

(ii) If both $abi, acj\in E(G)$, then $abi, acj$ and a $P_{t-2}$ in $[2t-2]\setminus \{k\}$ with endpoints $\{i, j\}$, where $k\in[3]\setminus \{i, j\}$, form a $C_t^3.$

(iii) If both $aij, abk\in E(G)$, then $aij, abk$ and a $P_{t-2}$ in $[2t-2]$ with endpoints $\{i, k\}$ or endpoints $\{j, k\}$ form a $C_t^3$.
\q

\begin{fact}\label{1}\label{t1}
Let $G$ be a dense and $C_t^3$-free $3$-graph, and let $[2t-2]^3\setminus\{123\}\subseteq G[A]$.  Then each $a\in V(G)\setminus A$ belongs to at most 1 good pair to $A$.
\end{fact}
\noindent{\em  Proof.}
If $a\in V(G)\setminus A$ belongs to 2 good pairs $\{a, b_1\}$, $\{a, b_2\}$ to $A$, then $ab_11,ab_22 \in E(G)$. This contradicts to \ref{3t1}(ii).
\q

\begin{fact}\label{2}\label{t2}
Let $G$ be a good graph to $A$ with $V(G)=\cup_{i=1}^s\{a_i,b_i\}\cup A$, where $s=\frac{|V(G)\setminus A|}{2}$ and let $K_{2t-2}^{3-}\subseteq G[A]$.  Then for each pair $\{i,j\}$, $N(a_i,a_j)\subseteq \{b_i,b_j\}$ and $N(a_i,b_j)\subseteq \{b_i,a_j\}$.
\end{fact}
\noindent{\em  Proof.}
If $k\in N(a_i, a_j)\cap A$, then $a_ia_jk,a_ib_ik'\in E(G)$ for $k'\in A\setminus\{k\}$, a contradiction to Fact \ref{3t1}(ii). So $N(a_i,a_j)\cap A=\emptyset$. Suppose that $x\in N(a_i,a_j)\subseteq V(G)\setminus (A\cup \{b_i,b_j\})$. If $t\geq 4$, then $a_ia_jx, a_ib_i1\in E(G)$. But $a_ib_i1, a_ia_jx, a_jb_j2$ and a path $P_{t-3}$ connecting 1 and 2 form a $C_t^3,$ a contradiction. If t=3, then $a_ia_jx, a_ib_i1, a_jb_j1$ form a $C_3^3$, a contradiction. So $N(a_i,a_j)\subseteq\{b_i, b_j\}$. Similarly, $N(a_i,b_j)\subseteq\{b_i, a_j\}.$
\q

\begin{defi}
For  good pairs $\{a_1,b_1\}$, $\{a_2,b_2\}$ to $A$, we say $\{a_2,b_2\}\geq\{a_1,b_1\}$ if $a_2b_2a_1$, $a_2b_2b_1\in E(G)$.
\end{defi}

\begin{fact}\label{3}\label{t3}
Let $G$ be a good graph to $A$ with $V(G)=\cup_{i=1}^s\{a_i,b_i\}\cup A$, where $s=\frac{|V(G)\setminus A|}{2}$, and let $K_{2t-2}^{3-}\subseteq G[A]$. Then either $\{a_i,b_i\}\geq\{a_j,b_j\}$ or $\{a_j,b_j\}\geq\{a_i,b_i\}$.
\end{fact}
\noindent{\em  Proof.}
By Fact \ref{t2}, $N(a_i,a_j)\subseteq \{b_i,b_j\}$. Without loss of generality, let $a_ib_ia_j\in E(G)$. By Fact \ref{t2}, $N(a_i,b_j)\subseteq \{b_i,a_j\}$. If $a_ib_jb_i\in E(G)$, then $\{a_i,b_i\}\geq\{a_j,b_j\}$. Else  $a_ib_ja_j\in E(G)$, by Fact \ref{t2}, $N(b_i,b_j)\subseteq \{a_i,a_j\}$. If $a_ib_ib_j\in E(G)$, then $\{a_i,b_i\}\geq\{a_j,b_j\}$. Otherwise $\{a_i,b_i\}\leq\{a_j,b_j\}$.
\q

\begin{fact}\label{t3'}\label{4.9}
Let $G$ be a good graph to $A$ with $V(G)=\cup_{i=1}^s\{a_i,b_i\}\cup A$, where $s=\frac{|V(G)\setminus A|}{2}$. If $\lambda{(G[A])}\leq\lambda{(K_{2t-1}^3)}$, then $\lambda{(G)}\leq\lambda{(K_{2t-1}^3)}$, equality holds if and only if $G=G[A]$ and $\lambda{(G[A])}=\lambda{(K_{2t-1}^3)}$
\end{fact}

\noindent{\em  Proof.} Let  $O_s$ be the $3$-graph whose vertex set is $\{a_{1}, b_{1}, a_{2}, b_{2},\dots , a_{s}, b_{s}\}$, and edge set is $\{a_{i}b_{i}a_{j},$ $a_{i}b_{i}b_{j}|i\neq j, 1\leq i, j\leq s.\}$ By Fact \ref{t2} and the definition of a good $3$-graph, we have
 $$E(G)\subset E(G[A])\cup E(O_s)\cup \{a_ib_ia| a\in A, 1\leq i,j\leq s\}.$$  Let the sums of weights of $\cup_{i=1}^s\{a_i, b_i\}$ in an optimal  weighting is $a$. Then
 \begin{equation}\label{lambdaG}
 \lambda{(G)}\leq \lambda{(O_s)}a^3+\frac{a^2}{4}(1-a)+\lambda{(K_{2t-1}^3)}(1-a)^3. \ \
 \end{equation}
 Let's prove the following fact.
\begin{fact} \label{123}
$\lambda{(O_s)}=\frac{1}{16}$.
\end{fact}
\noindent{\em  Proof.}
By Fact \ref{symmetry}, we may assume that $a_{i}, b_{i}$ have the same weighting (say $w_i$) in an optimal weighting. Then $$\sum_{i=1}^sw_i=\frac{1}{2},$$
and  $$\lambda{(O_s)}=\sum_{i=1}^sw_i^2(1-2w_i).$$
By Fact \ref{fact2}, we know that if $w_i, w_j>0$, then $$3\lambda{(O_s)}=w_i(1-2w_i)=w_j(1-2w_j),$$
i.e. $$(w_i-w_j)(1-2(w_i+w_j))=0.$$
So either $w_i=w_j$ or $w_i+w_j=\frac{1}{2}$. If $w_i+w_j=\frac{1}{2}$, then all other $w_k=0$  and this is $O_2$ which is isomorphic to $K_4^3$. It's easy to check that $\lambda{(O_2)}=\frac{1}{16}.$ So we only need to verify the case that $w_1=w_2=\dots=w_k=\frac{1}{2k}.$ In this case, $\lambda{(O_s, \overrightarrow{x})}=k(\frac{1}{2k})^2(1-\frac{1}{k})=\frac{k-1}{4k^2}\leq\frac{1}{16}$ for $k\geq 2.$
\q

 Let us continue the proof of Fact \ref{4.9}. Applying Fact \ref{123} to (\ref{lambdaG}), we have
\begin{equation}\label{lambdaG2}
\lambda{(G)}\leq\frac{1}{16}a^3+\frac{a^2}{4}(1-a)+m(1-a)^3 \ \ (m=\lambda{(K_{2t-1}^3)}=\frac{(2t-2)(2t-3)}{6(2t-1)^2}).
\end{equation}

\textbf{Case 1}. $t\ge 4$.

By (\ref{lambdaG2}),
\begin{eqnarray}
\lambda{(G)}
&\leq&\frac{1}{12}a^3+\frac{a^2}{4}(1-a)+m(1-a)^3\nonumber\\
&=&-(\frac{1}{6}+m)a^3+(3m+\frac{1}{4})a^2-3ma+m=f(a).\nonumber
\end{eqnarray}
Then  $f'(a)=-(\frac{1}{2}+3m)a^2+(6m+\frac{1}{2})a-3m$, and
the zeros of $f'(a)$ are $\frac{6m}{6m+1}$ and $1$. Since $f'(a)$ opens down, then $f(a)$ decrease in $[0, \frac{6m}{6m+1}]$ and increase in $[\frac{6m}{6m+1}, 1]$. So $$f(a)\leq Max\{f(0), f(1)\}=\{m, \frac{1}{12}\}=m$$
since $m\geq \frac{5}{49}.$ So $\lambda{(G)}\leq m=\lambda{(K_{2t-1}^3)}$. Equality holds if and only if $a=0$ and $\lambda{(G[A])}=\lambda{(K_{2t-1}^3)}.$

\textbf{Case 2}. $t=3$.

In this case, $m={2 \over 25}$. In view of (\ref{lambdaG2}),
\begin{eqnarray}
\lambda{(G)}
&\leq&\frac{1}{16}a^3+\frac{a^2}{4}(1-a)+{2 \over 25}(1-a)^3\nonumber\\
&=&\frac{1}{400}(-107a^3+196a^2-96a+32)=g(a)\nonumber.
\end{eqnarray}
Then  $g'(a)=\frac{1}{400}(-321a^2+392a-96)$ and
the zeros of $g'(a)$ are $x_{1}=(196-\sqrt{7600})/321$ and $x_{2}=(196+\sqrt{7600})/321$. Consequently, $g(a)$ decreases in $[0, x_1]$, increases in $[x_1, x_2]$ and decreases in $[x_2, 1]$. So $$g(a)\leq Max\{g(0), g(x_2)\}=\frac{2}{25}$$
and consequently $\lambda{(G)}\le \frac{2}{25}$. Equality holds if and only if $a=0$ and $\lambda{(G[A])}=\lambda{(K_5^3)}.$
\q

\bigskip
\subsection{Lagrangian density of $C_3^3$}
 We give the proof of Theorem \ref{c33}  in this section. Let $G$ be a dense and $C_3^3$-free $3$-graph, we show that $G$ is a good graph satisfying Fact \ref{4.9} and obtain the conclusion by applying Fact \ref{4.9}.

\noindent{\em  Proof of Theorem \ref{c33}.}
By Remark \ref{remark}, we may assume that $G$ is dense. By a result of de Caen in \cite{caen}, $\pi(K_4^{3-})\le \frac{1}{3}$. Since $K_4^{3-}$ cover pairs, by Proposition \ref{prop1}, $\pi_{\lambda}(K_4^{3-})\le \frac{1}{3}.$ If $\lambda{(G)}>\frac{2}{25}$, then $G$ contains $K_4^{3-}.$ Let us assume that $\{123, 124, 134\}\subseteq G.$ We will show that $G$ is a good graph and apply Fact \ref{4.9}.

\textbf{Case 1}  Assume $G$ doesn't contain an isomorphic copy of $\{123,124,134,x12\}.$

\textbf{Case 1.1}  $G$ contains $\{123,124,134,x23\}.$

Then $xy1\in E(G)$ for $y\in V(G)\setminus \{1,2,3,4,x\}$, otherwise $xy1,x23,124$ form a $C_3^3$. So, $N(x,1)\subseteq \{2,3,4\}$, however this contradicts to the assumption of Case 1.

\textbf{Case 1.2}
 $G$ doesn't contain a copy isomorphic to $\{123,124,134,x23\}.$

In this case, we have $x23,x24,x34,x12,x13,x14\notin E(G)$ for any $x$. Then there exists $y\notin \{2,3,4\}$ such that $xy1\in E(G)$. We claim that $N(x,2)=N(x,3)=N(x,4)=\{y\}$. Otherwise, let $y_1\in N(x,2)$ and
$y\neq y_1$, clearly $y_1\notin \{1, 3, 4\}$. Then $xy_12,xy1,123$ form a $C_3^3$, a contradiction. Similarly, the same holds for $N(x,3)$ and $N(x,4)$. Let $A=\{1,2,3,4\}$, then what we have obtained is that for each $x$ there exists $y$ such that $\{x,y\}$ is a good pair to $A$, so $G$ is a good graph to $A$.  Since $\lambda{(G[\{1,2,3,4\}])}<\frac{2}{25}$, then $\lambda{(G)}<\frac{2}{25}$ by Fact \ref{4.9}.

\textbf{Case 2}
 Assume that $G$ contains $\{x12,123,124,134\}$ for some $x\in V(G).$

\textbf{Case 2.1}
  $x23\in E(G)$.

In this case we have the following observation.
\begin{obser} \label{obser4.11}
For any $a\in V(G)\setminus \{1,2,3,4,x\}$, there exists  $b\in V(G)\setminus \{1,2,3,4,x\}$ such that $ab3\in E(G)$.
\end{obser}
\noindent{\em  Proof.}
Note that $a13\notin E(G)$ since otherwise $a13,x23,124$ form a $C_3^3$,  $a23\notin E(G)$ since otherwise $a23,x12,134$ form a $C_3^3$,  $a34\notin E(G)$ since otherwise $a34,x23,124$ form a $C_3^3$, and $ax3\notin E(G)$ since otherwise $ax3,x12,134$ form a $C_3^3$. So there exists  $b\in V(G)\setminus\{1, 2, 3, 4, x\}$ such that $b\in N(a,3)$ since $G$ covers pair.
\q

Now we fix such a $b$.

\begin{obser}\label{o4.12}
$N(a,1)=N(a,2)=N(a,3)=N(a,4)=N(a,x)=\{b\}$.
\end{obser}
\noindent{\em  Proof.}
Note that $a12\notin E(G)$  since otherwise $a12,ab3,134$ form a $C_3^3$,  $a13\notin E(G)$  since otherwise $a13,x23,124$ form a $C_3^3$,  $a14\notin E(G)$  since otherwise $a14,ab3,123$ form a $C_3^3$,  $ax1\notin E(G)$  since otherwise $ax1,ab3,134$ form a $C_3^3$, and $ay1\notin E(G)$  for $y\in V(G)\setminus \{1,2,3,4,x,b\}$ since otherwise $ay1,ab3,134$ form a $C_3^3$. So  $N(a,1)=\{b\}$.

We have shown that $a12, a23\notin E(G)$. Note that  $a24\notin E(G)$  since otherwise $a24,ab3,123$ form a $C_3^3$,  $ax2\notin E(G)$  since otherwise $ax2,ab3,123$ form a $C_3^3$, and $ay2\notin E(G)$  for $y\in V(G)\setminus \{1,2,3,4,x,b\}$ since otherwise $ay2,ab3,123$ form a $C_3^3$. So $N(a,2)=\{b\}$.

We have shown that $a14, a24, a34\notin E(G)$. Note that  $ax4\notin E(G)$  since otherwise $ax4,ab3,134$ form a $C_3^3$, and $ay4\notin E(G)$  for $y\in V(G)\setminus \{1,2,3,4,x,b\}$ since otherwise $ay4,ab3,134$ form a $C_3^3$. So  $N(a,4)=\{b\}$.

We have shown that $\{ax1, ax2, ax3, ax4\}\cap E(G)=\emptyset$. Note that  $axy\notin E(G)$ for $y\in V(G)\setminus \{1,2,3,4,x,b\}$ since otherwise $axy,ab3,x23$ form a $C_3^3$. So $N(a,x)=\{b\}$.

By Observation \ref{obser4.11} we know that $b\in N(a, 3)$. Assume that there exist $b'\neq b$  such that $ab_13\in E(G)$. Then $ab_13, ab1, 123$ form a $C_3^3.$ So $N(a, 3)=\{b\}.$
\q

By Observation \ref{o4.12}, $\{a,b\}$ is a good pair to $\{1,2,3,4,x\}$,  therefore $G$ is a good graph to $\{1, 2, 3, 4, x\}$. Since $\lambda{(G[\{1,2,3,4,x\}])}\leq\frac{2}{25}$, then $\lambda{(G)}\leq\frac{2}{25}$ by Fact \ref{4.9}. Equality holds if and only if $G[\{1, 2, 3, 4, x\}]=K_5^3.$

\textbf{Case 2.2} $G$ doesn't contain a copy isomorphic to $\{x23,x12,123,134,124\}.$

In this case,we know that $x23,x24\notin E(G)$. Then $N(x,3)\subseteq\{1,4\}$ since otherwise let $y\in N(x,3)$ and $y\in V(G)\setminus \{1,2,3,4,x\}$, then $xy3,x12,134$ form a $C_3^3$.

\textbf{Case 2.2.1} $x13\in E(G)$.

\textbf{Case 2.2.1.1} There exists $a\notin \{1,2,3,4,x\}$ such that $ax1\in E(G)$.

In this case, we know that $ax2\notin E(G)$  since otherwise $\{x12,123,x13,ax1,ax2\}$ is isomorphic to $\{123,124,134,x12,x23\}$. And $a12,a13,a14\in E(G)$  since otherwise if $1\in N(a,2)$, then there exist $y\notin \{1,2,x\}$ such that $ay2\in E(G)$, however $ax1,ay2,123$ form a $C_3^3$ (change 123 to 124 if $y=3$). Similarly the same holds for $N(a,3)$ and $N(a,4)$. We claim that $x14\in E(G)$  since otherwise there exists $z\in N(x,4)\setminus \{1\}$, then $xz4,ax1,134$ form a $C_3^3$ (change 134 to 124 if $z=3$). Let $A$ be a maximal set containing $\{2, 3, 4, x, a\}$ such that $1A^2\subseteq E(G)$ i.e. $1ij\in E(G)$ for all $i,j\in A$ and $i\neq j$. We have the following claim.
\begin{claim}
For any $a_i\in V(G)\setminus (A\cup\{1\})$, there  exists exactly one $b_i$ such that $\{a_i,b_i\}$ is a good pair to $A\cup\{1\}$ of $G$.
\end{claim}
\noindent{\em  Proof.}
 We claim that $a_i1k\notin E(G)$ for any $k\in A$ and $a_i\in V(G)\setminus (A\cup\{1\})$. If the conclusion does not hold i.e. there exist $a_i\in V(G)\setminus (A\cup\{1\})$ such that $a_i1j\in E(G)$ for some $j\in A$.  Without loss of generality, let $a_i12\in E(G)$. Consider $N(a_i, 3)$, note that $t\notin N(a_i, 3)$ for $t\in A$ since otherwise $a_it3, 1ta, 13x$ form a $C_3^3$ (change $a$ to 2 if $t=a$ or change $x$ to 2 if $t=x$).  And $N(a_i, 3)\cap(V(G)\setminus (A\cup\{1\}))=\emptyset$ since otherwise if $y\in N(a_i, 3)\cap(V(G)\setminus (A\cup\{1\}))$ then $a_iy3, a_i12, 134$ form a $C_3^3$. So $N(a_i, 3)=\{1\}$.  Similarly $N(a_i, k)=\{1\}$ for all $k\in A$, a contradiction to the maximality of $A$. So for any $a_i\in V(G)\setminus (A\cup\{1\})$ there exist $b_i\in V(G)\setminus (A\cup\{1\})$ such that $b_i\in N(a_i, 1)$. And $b_i\in N(a_i, k)$ for $k\in A$ since otherwise if $b_i\notin N(a_i, k)$ for some $k$, then there exists $b_i'\in N(a_i, k)$, consequently $a_ib_i'3, a_ib_i1, 134$ form a $C_3^3$. So $N(a_i, k)=\{b_i\}$. If there exist $b_i''\neq b_i$ such that $b_i''\in N(a_i, 1)$, then $a_ib_i''1, a_ib_i2, 123$ form a $C_3^3$, a contradiction. So $N(a_i, 1)=b_i$ as well and we have shown that $\{a_i, b_i\}$ is a good pair to $A\cup\{1\}.$
\q

By the definition of $A$,  $S_k^3(1)\subset G[A\cup \{1\}]$, where $k=\vert A\vert +1$. And any edge in $A$ would lead to a $C_3^3$, so $G[A\cup \{1\}]=S_k^3(1)$ and it's easy to get that $\lambda{(G[A\cup \{1\}])}<\frac{2}{27}$.
By Fact \ref{4.9}, $\lambda{(G)}<\frac{2}{25}$.

\textbf{Case 2.2.1.2} For any $a\in V(G)\setminus \{1,2,3,4,x\}$, $ax1\notin E(G)$.

In this case we claim that
\begin{claim}\label{4.14}
There exists $b\in V(G)\setminus \{1,2,3,4,x\}$ such that $\{a,b\}$ is a good pair to $\{1,2,3,4,x\}$.
\end{claim}
\noindent{\em  Proof.}
 Let's consider $N(a,x)$. Note that $2\notin N(a,x)$ since otherwise $ax2,x13,124$ form a $C_3^3$, $3\notin N(a,x)$ since otherwise $ax3,x12,134$ form a $C_3^3$, and $4\notin N(a,x)$ since otherwise $ax4,x12,134$ form a $C_3^3$. So there exists $b\in V(G)\setminus \{1,2,3,4,x\}$ such that $abx\in E(G)$.

Consider $N(a,1)$. Note that $2\notin N(a,1)$ since otherwise $abx,a12,x13$ form a $C_3^3$,  $3\notin N(a,1)$  since otherwise $abx,a13,x12$ form a $C_3^3$,  $4\notin N(a,1)$ since otherwise $abx,a14,x12$ form a $C_3^3$, and $N(a,1)\cap(V(G)\setminus \{1,2,3,4,x,b\})=\emptyset$   since otherwise $abx,ay1,x12$ form a $C_3^3$ for some $y\in V(G)\setminus \{1,2,3,4,x,b\}$. So $N(a,1)=\{b\}$.

If $N(a,2)\neq \{b\}$, then there exists $b'\in N(a,2)$. Note that $b'\neq 1,x$ since we have shown that $a12,ax2\notin E(G)$, then $ab'2,abx,x12$ form a $C_3^3$. So $N(a,2)=\{b\}$.

If $N(a,3)\neq \{b\}$, then there exists $b''\in N(a,3)$. Note that $b''\neq 1,x$ since we have shown that $a13,ax3\notin E(G)$, then $ab''3,abx,x13$ form a $C_3^3$. So $N(a,3)=\{b\}$.

If $N(a,4)\neq \{b\}$, then there exists $b'''\in N(a,4)$. Note that $b'''\neq 1,3$ since we have shown that $a14,a34\notin E(G)$, then $ab'''4,ab1,134$ form a $C_3^3$. So $N(a,4)=\{b\}$.

Therefore $\{a, b\}$ is a good pair to $\{1, 2, 3, 4, x\}.$
\q

By Claim \ref{4.14}, $G$ is a good graph to $\{1, 2, 3, 4, x\}$. Since $\lambda{(G[1,2,3,4,x])}\leq\frac{2}{25}$, then  $\lambda{(G)}\leq\frac{2}{25}$  by Fact \ref{4.9}. Equality holds if and only if $G[\{1, 2, 3, 4, x\}]=K_5^3.$

\textbf{Case 2.2.2} $x13\notin E(G)$.

Recall that $N(3,x)\subseteq \{1,4\}$, then $x34\in E(G)$.
\begin{obser}
For every $a\in V(G)\setminus \{1,2,3,4,x\}$ there exists $b\in V(G)\setminus \{1,2,3,4,x,a\}$ such that $N(a,x)\supseteq\{b\}$.
\end{obser}
\noindent{\em  Proof.}
Note that $ax1\notin E(G)$  since otherwise $ax1,x34,123$ form a $C_3^3$,  $ax2\notin E(G)$  since otherwise $ax2,x34,123$ form a $C_3^3$,  $ax3\notin E(G)$  since otherwise $ax3,x12,134$ form a $C_3^3$, and $ax4\notin E(G)$  since otherwise $ax4,x12,134$ form a $C_3^3$. So there exists $b\in V(G)\setminus \{1,2,3,4,x,a\}$ such that $abx\in E(G)$.
\q

\begin{claim}
Let $A=\{a:a12,a34\in E(G)\}$, then $G$ is a good graph to $A\cup \{1,2,3,4\}$.
\end{claim}
\noindent{\em  Proof.}
Note that $A\neq\emptyset$ since $x\in A$. We will prove an observation first.
\begin{obser}\label{obser4.17}
For any $y\in V(G)\setminus (A\cup\{1,2,3,4\})$, $y12\notin E(G)$.
\end{obser}
\noindent{\em  Proof.}
If $y12\in E(G)$, then consider $N(y,3)$. Take $a\in A$. If $1\in N(y,3)$, then $a12,y13,a34$ form a $C_3^3$. If $2\in N(y,3)$, then $y23,a12,134$ form a $C_3^3$. If $a\in N(y,3)$, then $ay3,a12,134$ form a $C_3^3$. If $x'\in N(y,3)$, where $x'\notin A\cup\{1, 2, 3, 4\}$, then $x'y3,y12,134$ form a $C_3^3$. So $N(y,3)=\{4\}$, then  $y\in A$,  a contradiction to the maximality of $A$.
\q

For any $z\in V(G)\setminus (A\cup\{1,2,3,4\})$, we consider $N(z,1)$, note that $2\notin N(z,1)$ by Observation \ref{obser4.17}. Take $a\in A$. Note that $3\notin N(z,1)$ since otherwise $z13,a34,a12$ form a $C_3^3$,  $4\notin N(z,1)$ since otherwise $z14,a34,a12$ form a $C_3^3$, and $a\notin N(z,1)$ since otherwise $az1,a34,123$ form a $C_3^3$. So there exists $u\in V(G)\setminus A\cup \{1,2,3,4\}$ such that $uz1\in E(G)$.

Consider $N(z,2)$. If $u\notin N(z,2)$, then there exists $v\in N(z,2)$ and $vz2,uz1, 123$ form a $C_3^3$ (change 123 to 124 if $v=3$). So $N(z,2)=\{u\}$. Similarly, we have $N(z,4)=\{u\}$, and $N(z,a)=\{u\}$ for $a\in A$. So $\{z,u\}$ is a good pair to $A\cup\{1,2,3,4\}$. Hence $G$ is a good graph to $A\cup\{1,2,3,4\}$.
\q

Since $G$ is a good graph to $A\cup\{1, 2, 3, 4\}$, then for $a\in A$, $N(a,1)\subseteq\{2, 3, 4\}$. If $a13\in E(G)$ or $a14\in E(G)$, then it is Case 2.2.1 since $a12\in E(G)$. So we may assume that $N(a,1)=\{2\}$. Similarly, we can show that $N(a,2)=\{1\}$. Consequently, $N(a,3)=\{4\}$ and $N(a,4)=\{3\}$. Hence $G$ is a good graph to $A$. Applying the same procedure to $A$. If $\lambda{(G[A])}<\frac{2}{25},$ then by Fact \ref{4.9}, we know that Theorem \ref{c33} holds. Else if $\lambda{(G[A])}\geq\frac{2}{25}$, then applying the same procedure to $G[A]$ implies that there exist 4 vertices $\{5, 6, 7, 8\}$ such that $\{5, 6\}$ and $\{7, 8\}$ are good pairs to $A^{(1)}=A\setminus\{5, 6, 7, 8\}$. Continue this procedure, until we obtain that $\lambda{(A^{(i)})}<\frac{2}{25}$ or $|A^{(i)}|\leq5$ and $G$ is a good graph to $A^{(i)}$. By Fact \ref{4.9}, $\lambda{(G)}\leq\frac{2}{25},$ and equality holds if and only if $G[\{1, 2, 3, 4, x\}]=K_5^3.$
\q

\bigskip
\subsection{Extension to $C_t^3$}

In this section we prove Theorem \ref{3t'}.
Let $G$ be a dense and $C_t^{3}$-free $3$-graph with $K_{2t-2}^{3-}\subset G$ $(t\geq 4)$. Without loss of generality, let $[2t-2]^3\setminus \{123\}$ be a $K_{2t-2}^{3-}$ in $G$. We will show that $G$ is a good graph satisfying Fact \ref{4.9} and obtain the conclusion of Theorem \ref{3t'} by applying Fact \ref{4.9}.

\begin{fact} \label{3t2}
Let $G$ be a dense and $C_t^{3}$-free $3$-graph. Let $[2t-2]^3\setminus \{123\}$ be a $K_{2t-2}^{3-}$ in $G$. If there exist $a, b\in V(G)\setminus [2t-2]$ such that $abi\in E(G)$ for some $i\in [2t-2]$, then $\{a, b\}$ is a good pair to $[2t-2]$ or $N(a,k)=N(b,k)=i$ for $k\in [2t-2]\setminus\{i\}.$
\end{fact}
\noindent{\em  Proof.}
For $k\in [2t-2]\setminus\{i\}$, if $N(a,k)\neq\{i\}$ and $N(a,k)\neq\{b\}$, then there exists $b'\notin\{i, b\}$ such that $ab'k\in E(G)$. But $abi, ab'k\in E(G)$ contradicting to Fact \ref{3t1}(ii) or (iii).

By Fact \ref{3t1}(iii), it's impossible that $N(a,k)=\{i\}$ and $N(a,k')=\{b\}$ for some $k, k'\in [2t-2]\setminus\{i\}.$ So either (1) $N(a,k)=\{i\}$ for all $k\in [2t-2]\setminus\{i\}$ or (2) $N(a,k)=\{b\}$ for all $k\in [2t-2]\setminus\{i\}$. So what left is to show that $N(a,i)=\{b\}$ if $N(a,k)=\{b\}$ for all $k\in [2t-2]\setminus\{i\}.$ For $k\in [2t-2]\setminus\{i\}$, by Fact \ref{3t1}(iii), $aki\notin E(G)$ since $abk'\in E(G)$ for $k'\in [2t-2]\setminus\{i, k\}$. For $c\in V(G)\setminus([2t-2]\cup\{a, b\})$, by Fact \ref{3t1}(ii), $aci\notin E(G)$ since $abk\in E(G)$. So $N(a,i)=\{b\}.$
\q

\noindent{\em  Proof of Theorem \ref{3t'}.}
Let $G$ be a dense and $C_t^{3}$-free $3$-graph with $K_{2t-2}^{3-}\subset G$ $(t\geq 4)$. If $|V(G)|\leq2t-1$, then $\lambda{(G)}\leq\lambda{(K_{2t-1}^3)}$, equality holds if and only if $G=K_{2t-1}^3.$ Assume that $|V(G)|\geq2t$.  Without loss of generality, let $[2t-2]^3\setminus \{123\}$ be a $K_{2t-2}^{3-}$ in $G$. We will show that $G$ is a good graph satisfying Fact \ref{4.9} and obtain the conclusion of Theorem \ref{3t'} by applying this fact.

\textbf{Case 1} $G$ doesn't contain a copy isomorphic to $\{x14, [2t-2]^3\setminus \{123\}\},$ where $x\in V(G)\setminus [2t-2].$

In this case, we know that $xk1,xk2,xk3\notin E(G)$ for $4\leq k\leq 2t-2$ and $x\in V(G)\setminus [2t-2].$

\textbf{Case 1.1} $G$ contains a copy isomorphic to $\{x45, [2t-2]^3\setminus \{123\}\}.$

\textbf{Case 1.1.1} $G$ doesn't contain a copy isomorphic to $\{x45, x46, [2t-2]^3\setminus \{123\}\}.$

Assume that $x45\in E(G)$, then $x46\notin E(G)$. If $y\in N(x,6)\cap (V(G)\setminus [2t-2])$, then $xy6, x45, 612$ together with a $P_{t-3}$ in $[2t-2]\setminus \{2, 5, 6\}$ connecting 1 and 4 form a $C_t^3$. Without loss of generality, let $x67\in E(G)$, then similarly we have $x89, x(10)(11), ... , x(2t-4)(2t-3)\in E(G),$ however $N(x, 2t-2)=\emptyset$, a contradiction.

\textbf{Case 1.1.2} $G$ contains a copy isomorphic to $\{x45, x46, [2t-2]^3\setminus \{123\}\}.$

Without loss of generality, assume that $x45, x46\in E(G)$. Then for any $a\in V(G)\setminus [2t-2]$, $x\notin N(a, 5)$  since otherwise $xa5, x46$ together with a $P_{t-2}$ on $[2t-2]\setminus \{4\}$ connecting 5 and 6 form a $C_t^3$. Moreover $k\notin N(a, 5)$ for $k\in [2t-2]\setminus \{5\}$, otherwise $ak5, x45$ together with a $P_{t-2}$ on $[2t-2]\setminus \{5\}$ connecting 4 and $k$ form a $C_t^3$ (change $x45$ to $x46$ if $k=4$). Then there exist $b\in V(G)\setminus V(K_{2t-2}^{3-})$ such that $ab5\in E(G)$. By Fact $\ref{3t2}$, $\{a, b\}$ is a good pair to $[2t-2]$ since $N(a,k)\neq\{5\}$ for $i\in [2t-2]$.  Then By Fact \ref{t3'}, $\lambda{(G)}<\lambda{(K_{2t-1}^3)}$.

\textbf{Case 1.2} $G$ doesn't contain a copy isomorphic to $\{x45, [2t-2]^3\setminus \{123\}\}.$\\
In this case, for any $a\in V(G)\setminus [2t-2]$, there exists $b\in V(G)\setminus [2t-2]$ such that $abi\in E(G)$ and $N(a,k)\neq\{i\}$ for any $k\in [2t-2]\setminus\{i\}$. By Fact \ref{3t2}, $\{a, b\}$ is a good pair to $[2t-2]$. By Fact \ref{t3'}, $\lambda{(G)}< \lambda{(K_{2t-1}^3)}$.

\textbf{Case 2} $G$ contains a copy isomorphic to $\{x41, [2t-2]^3\setminus \{123\}\}.$

\textbf{Case 2.1} $G$ contains a copy isomorphic to $\{x41, x45, [2t-2]^3\setminus \{123\}\}.$

Without loss of generality, assume that $x41, x45\in E(G)$. Then for any $a\in V(G)\setminus [2t-2]$, we know that $x\notin N(a, 5)$ since otherwise $x41, xa5$ together with a $P_{t-2}$ on $[2t-2]\setminus \{1\}$ connecting 4 and 5 form a $C_t^3$. Moreover $k\notin N(a, 5)$ for $k\in [2t-2]\setminus \{5\}$, since otherwise $ak5, x45$ together with a $P_{t-2}$ on $[2t-2]\setminus \{5\}$ connecting 4 and $k$ form a $C_t^3$ (change $x45$ to $x41$ if $k=4$). Then there exist $b\in V(G)\setminus V(K_{2t-2}^{3-})$ such that $ab5\in E(G)$. Since we have already known that $N(a, k)\neq\{5\}$, then by Fact $\ref{3t2}$, $\{a, b\}$ is a good pair to $[2t-2]$. Applying Fact \ref{t3'}, we have $\lambda{(G)}<\lambda{(K_{2t-1}^3)}$.

\textbf{Case 2.2} $G$ contains a copy isomorphic to $\{x41, x15, [2t-2]^3\setminus \{123\}\}.$\\
Without loss of generality, we can assume that $x15, x41\in E(G).$ Then for any $a\in V(G)\setminus [2t-2]$, we know that $x\notin N(a, 5)$ since otherwise $x41, xa5$ together with a $P_{t-2}$ on $[2t-2]\setminus \{1\}$ connecting 4 and 5 form a $C_t^3$. Moreover $k\notin N(a, 5)$ for $k\in [2t-2]\setminus \{5\}$, since otherwise $ak5, x15$ together with a $P_{t-2}$ on $[2t-2]\setminus \{5\}$ connecting 1 and $k$ form a $C_t^3$ (change $x15$ to $x41$ if $k=1$, and find the linear path $P_{t-2}$ on $[2t-2]\setminus \{1\}$ connecting 4 and $5$). Then there exist $b\in V(G)\setminus V(K_{2t-2}^{3-})$ such that $ab5\in E(G)$. Since we have already known that $N(a, k)\neq\{5\}$, then by Fact $\ref{3t2}$, $\{a, b\}$ is a good pair to $[2t-2]$. Applying Fact \ref{t3'}, we have $\lambda{(G)}<\lambda{(K_{2t-1}^3)}$.

\textbf{Case 2.3} $G$ does not contain a copy isomorphic to $\{x41, x45, [2t-2]^3\setminus \{123\}\}.$

Without loss of generality, we assume that $x41\in E(G)$, $x15, x45\notin E(G)$. Note that $N(x, 5)\cap (V(G)\setminus [2t-2])=\emptyset$. Otherwise let $y\in N(x, 5)\cap (V(G)\setminus [2t-2])$, then $xy5, x41$ together with a $P_{t-2}$ on $[2t-2]\setminus \{1\}$ connecting 4 and 5 form a $C_t^3$. So there exist $k\in [2t-2]\setminus \{1, 4\}$ such that $xk5\in E(G)$. Then for any $a\in V(G)\setminus [2t-2]$, we have $x, 1, 4\notin N(a, 5)$ since otherwise $x41, ai5$ ($i=x, 1, 4$) together with a $P_{t-2}$ in $[2t-2]$ form a $C_t^{3}$. Note that $k\neq1$, then $2, 3, k\notin N(a, 5)$ since $x41, ai5\in E(G)$ ($i=2, 3, k$), contradicting to Fact \ref{3t1}(i). Moreover $y\notin N(a, 5)$ for $y\in [2t-2]\setminus \{1, 2, 3, 4, 5, k\}$ since otherwise $ay5, x5k$ together with a $P_{t-2}$ in $[2t-2]\setminus \{5\}$ form a $C_t^3$. So there exists $b\in V(G)\setminus [2t-2]$ such that $ab5\in E(G)$. Since $5\notin N(a, k)$, then by Fact \ref{3t2}, $\{a, b\}$ is a good pair to $[2t-2]$. By Fact \ref{t3'},  $\lambda{(G)}<\lambda{(K_{2t-1}^3)}$.
\q

\bigskip




\begin{thebibliography}{JluR00}

\bibitem{NWY}  A. Bene Watts, S. Norin and L. Yepremyan,  A Tur\'an theorem for extensions via an Erd\H {o}s-Ko-Rado theorem for Lagrangians, arXiv:1707.01533.

\bibitem{BIJ}
A. Brandt, D. Irwin, and T. Jiang, Stability and Tur\'an numbers of a class of hypergraphs via
Lagrangians,  \emph{Combin., Probab. \& Comput.},  \textbf{26(3)} (2017), 367-405.

\bibitem{CLP}
P. Chen, J. Liang and Y. Peng, The Lagrangian density of the disjoint union of a 3-uniform tight path and a matching and the Tur\'an number of its extension, manuscript.

\bibitem{caen} D. de Caen, Extension of a theorem of Moon and Moser on complete subgraph, \emph{Ars Combinatoria} 16 (1983) 5-10.





\bibitem{FF}
P. Frankl, Z. F\"uredi, Extremal problems whose solutions are the
blow-ups of the small Witt-designs, {\it J. Combin. Theory Ser. A},
{\bf 52}(1989), 129--147.


\bibitem{FR}
P. Frankl, V. R\"odl, Hypergraphs do not jump, {\it Combinatorica},
{\bf 4} (1984), 149--159.

\bibitem{GLM2018}V. Gruslys, S. Letzter, N. Morrison, Hypergraph Lagrangians: Resolving the Frankl-F\"uredi
conjecture, arXiv:1807.00793v1, preprint, JulY 2018.

\bibitem{HK}
D. Hefetz and P. Keevash, A hypergraph Tur\'an theorem via Lagrangians of intersecting
families, {\it J. Combin. Theory Ser. A}, {\bf 120} (2013), 2020--2038.

\bibitem{HPW}
S. Hu, Y. Peng and B. Wu, Lagrangian densities of linear forests and Tur\'an numbers of their extensions, submitted.



\bibitem{Jenssen} M. Jenssen, Continous Optimisation in Extremal Combinatorics, Ph.D. dissertation, London School of Economics and Political Science, 2017.

    \bibitem{JPW} T. Jiang, Y. Peng, B. Wu, Lagrangian densities of some sparse hypergraphs and Tur\'an numbers of their extensions, {\it European Journal of Combinatorics} {\bf 73}(2018), 20-36.

\bibitem{KNS}
G. Katona, T. Nemetz and M. Simonovits, On a problem of Tur\'an in the theory of graphs, \emph {Mat. Lapok}, {\bf 15} (1964), 228--238.

\bibitem{Keevash} P. Keevash, Hypergrah Tur\'an problems, {\it Surveys in Combinatorics}, Cambridge
University Press, {\bf392} (2011), 83--140.


\bibitem{LLP2018}H. Lei, L. Lu, and Y. Peng, On Lagrangians of 3-uniform hypergraphs,
arXiv:1806.10846v1, preprint, June 2018.

\bibitem{MS} T.S. Motzkin and E.G. Straus, Maxima for graphs and a new proof of a theorem of Tur\'an, \emph{Canad. J. Math}, {\bf 17} (1965), 533--540.




\bibitem{NY} S. Norin, L. Yepremyan, Tur\'an numbers of generalized triangles, \emph{J. Combin. Theory Ser. A}, 146 (2017), 312-343.

\bibitem{NY2} S. Norin, L. Yepremyan, Tur\'an numbers of extensions, \emph{J. Combin. Theory Ser. A}, 155 (2018), 476-492..

\bibitem{Pikhurko} O.Pikhurko, An exact Tur\'an result for the generalized triangle, \emph{Combinatorica} {\bf 28} (2008), 187--208.

\bibitem{Sidorenko-87} A.F. Sidorenko, On the maximal number of edges in a homogeneous hypergraph that does not contain prohibited subgraphs,
\emph{Mat. Zametki}, \textbf{41} (1987), 433-455.

\bibitem{Sidorenko-89} A.F. Sidorenko, Asymptotic solution for a new class of forbidden $r$-graphs, \emph{Combinatorica}, {\bf 9} (1989), 207--215.

\bibitem{T} J. Talbot, Lagrangians of hypergraphs, \emph{Combin., Probab. \& Comput.}, {\bf 11} (2002), 199--216.

\bibitem{TPZZ2} Q. S. Tang, Y. Peng, X. D. Zhang, and C. Zhao, Connection between the clique number and the Lagrangian of $3$-uniform hypergraphs, \emph{Optimization Letters}, 10.4 (2016), 685-697.

\bibitem{Tyo} M. Tyomkyn, Lagrangians of hypergraphs: The Frankl-F\"uredi conjecture holds almost everywhere, \emph{J.London Math. Soc.} {\bf 96} (2017), 584-600.



\bibitem{WPC} B. Wu, Y. Peng, and P. Chen, On a conjecture of Hefetz and Keevash on the Lagrangian density of intersecting hypergraphs, arXiv:1701.06126.

\bibitem{WP} B. Wu and Y. Peng, Lagrangian densities of short 3-uniform linear paths and Tur\'an number of their extensions, manuscript.




\bibitem{Z} A. A. Zykov, On some properties of linear complexes, \emph{Mat. Sbornik, (N. S.)}, \textbf{24} )(66) (1949), 163--188.


\end{thebibliography}
\end{document}